\documentclass[11pt]{amsart}

\usepackage{times}
\usepackage{amsmath,  amssymb,slashed,url,bm,upgreek,amsthm}
\usepackage{graphicx,enumerate}
\newcounter{intro}

\newtheorem{theo}[intro]{Theorem}
\newtheorem*{theo*}{Theorem}
\newtheorem{coro}[intro]{Corollary}

\newtheorem{thm}{Theorem}[section]
\newtheorem{lem}[thm]{Lemma}

\newtheorem*{merci}{Acknowledgements}

\newcommand{\cref}[1]{Corollary~\ref{#1}}

\newcommand{\lref}[1]{Lemma~\ref{#1}}

\newcommand{\tref}[1]{Theorem~\ref{#1}}

\DeclareMathOperator{\ricci}{Ricci}

\DeclareMathOperator{\diam}{diam}

\DeclareMathOperator{\vol}{dvol}
\DeclareMathOperator{\volu}{vol}
\def\R{\mathbb R}

\def\cC{\mathcal C}

\def\cH{\mathcal H}

\def\cO{\mathcal O}

\begin{document}
\title{Harmonic functions on Manifolds whose large sphere are small .}
\author{ Gilles Carron }

\address{Laboratoire de Math\'ematiques Jean Leray (UMR 6629), Universit\'e de Nantes, 
2, rue de la Houssini\`ere, B.P.~92208, 44322 Nantes Cedex~3, France}
\email{Gilles.Carron@univ-nantes.fr}
\begin{abstract}
 
 We study the growth of harmonic functions on complete Riemannian manifolds  where the extrinsic diameter of geodesic spheres is  sublinear. It is an generalization of a result of A. Kazue. We also get a Cheng and Yau estimates for the gradient of harmonic functions.
\\
R\'ESUM\'E :  On \'{e}tudie la croissance des  fonctions harmoniques sur les vari\'{e}t\'{e}s riemanniennes compl\`{e}tes dont le diam\`{e}tre des grandes sph\`{e}des g\'{e}od\'{e}siques  croit sous lin\'{e}airement. Il s'agit de g\'{e}n\'{e}ralisation de travaux de A. Kazue. Nous obtenons aussi une estim\'{e}e de type Cheng-Yau pour le gradient des fonctions harmoniques.
\end{abstract}
\maketitle
\section{Introduction}When $(M,g)$ is a complete Riemannian manifold with non negative Ricci curvature,
S-Y. Cheng and S-T. Yau have proven that any harmonic function $h\colon M\rightarrow \R$ satisfies the gradient estimate \cite{ChengYau} :
$$\sup_{z\in B(x,R)} |dh|(z)\le \frac{C(n)}{R}\sup_{z\in B(x,2R)} |h(z)|.$$
This result implies that such a manifold can not carry non constant harmonic function $h\colon M\rightarrow \R$ with sublinear growth :
$$|h(x)|=o\big(d(o,x)\big)\,\,,\,\, d(o,x)\to +\infty\,.$$
A celebrated conjecture of S-T. Yau predicted the finite dimensionality of the space of harmonic functions with 
polynomial growth on a complete Riemannian manifold with non negative Ricci curvature :
$$\cH_\nu(M,g)=\left\{h\in \cC^2(M)\,,\, \Delta_g h=0, |h(x)|=\cO\big(d^\nu(o,x)\big)\right\}.$$
This conjecture has been proven by T. Colding and B. Minicozzi in a much more general setting.

We say that  a complete Riemannian manifold $(M^n,g)$  satisfies the\textit{ doubling } condition if
 there is a constant $\upvartheta$ such that for any $x\in M$ and radius $R>0$ :
$$\hspace{1cm} \volu B(x,2R)\le  \upvartheta\volu B(x,R).$$
If $B\subset M$ is a geodesic ball, we will use the notation $r(B)$ for the radius of $B$ and  $\kappa B$ for 
the ball concentric to $B$ and with radius $\kappa r(B)$. And if $f$ is an integrable function on a subset 
$\Omega\subset M$, we will note $f_\Omega$ its mean over $\Omega$:
$$f_\Omega=\frac{1}{\volu \Omega} \int_\Omega f.$$

We say that  a complete Riemannian manifold $(M^n,g)$  satisfies  the scale ($L^2$) Poincar\'e inequality
 if there is a constant $\mu$ such that 
for any ball $B\subset M$ and any function $\varphi\in \cC^1(2B)$:
$$\|\varphi-\varphi_B\|^2_{L^2(B)}\le \mu\, r^2(B) \|d\varphi\|^2_{L^2(2B)}\,.$$
\begin{theo*}\cite{CM} If $(M,g)$ is a complete Riemannian manifold that is doubling and that satisfies  the scale
 Poincar\'e inequality then for any $\nu$, the space of harmonic function of polynomial growth of order $\nu$ has finite dimension:
$$\dim \cH_\nu(M,g)<+\infty.$$
\end{theo*}
It is well known that a complete Riemannian 
manifold with non negative Ricci curvature is doubling and satisfies the scale Poincar\'e inequality, hence the Yau's conjecture is true.

The proof is quantitative and gives a precise estimation of the dimension of $\dim \cH_\nu(M,g)$.
 In fact, the condition on the Poincar\'e inequality can be weakened and the result holds on a 
doubling manifold $(M,g)$  that satisfies the mean value estimation \cite{CM2,Li} : for any harmonic function defined over a geodesic ball $3B$ :
$$\sup_{x\in B} |h(x)|\le \frac{C}{\volu 2B} \int_{2B} |h|.$$

An example of Riemannian manifold satisfying the above condition are Riemannian manifold $(M,g)$ 
that outside a compact set $(M,g)$ is isometric to the warped product
$$( [1,\infty)\times \Sigma, (dr)^2+r^{2\gamma}h)$$ where $(\Sigma,h)$ is a closed connected manifold and 
$\gamma\in (0,1]$. But when $\gamma\in (0,1)$,  a direct analysis, separation of variables,
 shows that any harmonic function $h$ satisfying for some $\epsilon>0$: 
$$h(x)=\cO\left(e^{Cr^{1-\gamma-\epsilon}}\right)$$ is necessary constant. In particular, a harmonic 
function with polynomial growth is constant.  In \cite{Ka1,Ka2}, A. Kasue has shown that this was a general result for
 manifold whose Ricci curvature satisfies a quadratic decay lower bound and whose geodesic spheres have sublinear growth (see also \cite{sor} for a related results):
\begin{theo*}If $(M,g)$ is complete Riemannian manifold with a based point $o$ whose Ricci curvature 
satisfies a quadratic decay lower bound:
$$\hspace{1cm} \,\,\, \ricci\ge - \frac{\kappa^2}{d^2(o,x)} g\,\,, $$ and whose geodeosic sphere have
 sublinear growth:
$$\diam \partial B(o,R)=o(R)\,\,,\,\, R\to +\infty$$
then any harmonic function with polynomial growth is constant.
\end{theo*}

Following A. Grigor'yan and L. Saloff-Coste \cite{GS0}, we say that a ball $B(x,r)$ is remote (from a fixed point $o$) if
$$3r\le d(o,x).$$
Our first main result is a refinement of A. Kasue's result when the hypothesis of the Ricci curvature 
is replaced by a scale Poincar\'e inequality for remote ball :
There is a constant $\mu$ such that  all remote balls $B=B(x,r)$ satisfy a scale Poincar\'e inequality :
$$\forall \varphi\in \cC^1(2B)\,:\, \|\varphi-\varphi_B\|^2_{L^2(B)}\le \mu r^2 \|d\varphi\|^2_{L^2(2B)}$$

\begin{theo} Let $(M,g)$ be a complete Riemannian manifold whose remote balls satisfy the scale 
Poincar\'e inequality and assume  that
geodesic spheres have sublinear growth:
$$\diam \partial B(o,R)=o(R)\,\,,\,\, R\to +\infty.$$
If $h\colon M\rightarrow \R$ is a harmonic function such that for $I_R:=\int_{B(o,R)} h^2$ :
$$\liminf_{R\to +\infty}\,\, \log(I_R)\,\frac{\diam  \partial B(o,R)}{R}=0$$
then $h$ is constant.
\end{theo}
For instance, on such a manifold, a harmonic function $h\colon M\rightarrow \R$ satisfying :
$$|h(x)|\le C d(o,x)^\nu \left(\volu B(o,d(o,x))\right)^{-\frac12}$$ is  constant. Moreover if the diameter of geodesic sphere satisfies
$$\diam \partial B(o,R)\le C R^\gamma\,\,,$$
for some $\gamma\in (0,1)$ then 
if $h\colon M\rightarrow \R$ is a harmonic function such that for some positive constant $C$ and $\epsilon$:
$$|h(x)|\le C e^{d(o,x)^{1-\gamma-\epsilon}} \volu B(o,d(o,x))$$ then $h$ is  constant.

A by product of the proof will imply that on the class of manifold considered by A. Kasue, the
 doubling condition implies an estimate {\it\`a la }Cheng-Yau for  for the gradient of harmonic function:
\begin{theo}\label{lip}  Let $(M^n,g)$ be a complete Riemannian manifold that  is doubling
 and whose  Ricci curvature  satisfies a quadratic decay lower bound.
Assume that the diameter of geodesic sphere has a sublinear growth 
$$\diam \partial B(o,R)= \sup_{x,y\in \partial B(o,R)} d(x,y)=o(R)\,,$$

then there is a constant $C$ such that  for any geodesic ball $B\subset M$ and any harmonic function $h\colon 3B\rightarrow \R$
$$\sup_{x\in B} |dh|^2(x)\le \frac{C}{\volu 2B}\int_{2B} |dh|^2.$$
\end{theo}
This result has  consequences for the boundness of the Riesz transform. When $(M^n,g)$ is a
 complete Riemannian manifold with infinite volume, the Green formula and  the spectral theorem yield the equality:
$$\forall f\in \cC_0^\infty(M)\,\, ,\, \int_M |df|^2_g \vol_g=\langle \Delta f,f\rangle_{L^2}=
\int_M \left|\Delta^{\frac12}f\right|^2 \vol_g\,.$$
Hence the Riesz transform $$R:=d\Delta^{-\frac12}\colon L^2(M)\rightarrow L^2(T^*M)$$
 is a bounded operator. It is well known \cite{St2} that on a Euclidean space, 
the Riesz transform has a bounded extension $R\colon L^p(\R^n)\rightarrow L^p(T^*\R^n)$ for every $p\in (1,+\infty)$. 
Also according to D. Bakry, the same is true on manifolds with non-negative Ricci curvature \cite{Bakry}.
 As it was noticed in \cite[section 5]{CRiesz}, in the setting of the   \tref{lip}, 
the analysis of A. Grigor'yan and L. Saloff-Coste \cite{GS0} implies
a scale $L^1$-Poincar\'e inequality: 
there is a constant $C$ such that any balls $B=B(x,r)$ satisfies :
$$\forall \varphi\in \cC^1(2B)\,:\, \|\varphi-\varphi_B\|_{L^1(B)}\le C r^2 \|d\varphi\|_{L^1(2B)}\,. $$
And according to the analysis of P. Auscher and T. Coulhon \cite{AC} (see also the explanations in \cite[section 5]{CRiesz}), the  \tref{lip} implies :

\begin{coro}Under the assumption of \tref{lip},  the Riesz transform is bounded on $L^p$ for every $p\in (1,+\infty)$.
\end{coro}
\begin{merci} I thank Hans-Joachim Hein :  this project had begun by a  very  fruitful discussion where we proved together  the key lemma (\ref{inefunda}).
I'm partially supported by the grants ACG: ANR-10-BLAN 0105 and GTO : ANR-12-BS01-0004.
\end{merci}
\section{Absence of harmonic functions}

Recall that when $(M,g)$ is a complete Riemannian manifold and $o\in M$, we say that a geodesic ball $B(x,r)$  is {\it remote} (from $o$) if
$$3r\le d(o,x).$$
We define $\rho$ the radius function by 
$\displaystyle \rho(t)=\inf_{x\in \partial B(o,t)} \max_{y\in \partial B(o,t)} d(x,y)$, 
we have 
$$ \displaystyle\rho(t)\le \diam \partial B(o,t)\le 2\rho(t).$$

\subsection{An inequality}
\begin{lem} \label{inefunda}Let $(M,g)$ be a complete Riemannian manifold whose all remote balls $B=B(x,r)$ satisfy a scale Poincar\'e inequality :
$$\forall \varphi\in \cC^1(2B)\,:\, \|\varphi-\varphi_B\|^2_{L^2(B)}\le \mu\, r^2(B) \|d\varphi\|^2_{L^2(2B)}\, .$$
Then there are constants $C>0$ and $\kappa\in (0,1)$ depending only on $\mu$ such that if
$$\forall r\in [R,2R] : \rho(r)\le \varepsilon r \,\,\text{with}\,\, 
\varepsilon\in \left(0,1/12\right)$$  and if   $h$ is a harmonic function on $B(o,2R)$ then
$$\int_{B(o,R)} |dh|^2\le C\,\kappa^{\frac1\varepsilon}\, \int_{B(o,2R)} |dh|^2\,.$$\end{lem}
\proof Let $r\in [R+4\varepsilon R ,2R-4\varepsilon R]$, our hypothesis implies that there is some $x\in \partial B(o,r)$ such that 
$$B(o,r+\varepsilon R)\setminus B(o,r)\subset B(x,\varepsilon R+\varepsilon r).$$
Let $h\colon B(o,2R)\rightarrow \R$ be a harmonic function and $c\in \R$ a real number. We use the Lipschitz function :
$$\chi(x)=\begin{cases}
1&\text{on}\,\, B(o,r)\\
\frac{r+\varepsilon R-d(o,x)}{\varepsilon R}&\text{on}\,\, B(o,r+\varepsilon R)\setminus B(o,r)\\
0&\text{outside}\,\, B(o,r+\varepsilon R)\\
\end{cases}$$
Then integrating by part and using the fact that $h$ is harmonic we get 
$$\int_M \chi^2 |d(h-c)|^2+2\chi(h-c)\langle d\chi, d(h-c)\rangle =\int_M\langle d((h-c)\chi^2), d(h-c)\rangle=0$$
So that we have :\begin{equation*}\begin{split}\int_{M} |d(\chi(h-c))|^2&=\int_M \chi^2 |d(h-c)|^2+2\chi(h-c)\langle d\chi, d(h-c)\rangle+ (h-c)^2 |d\chi|^2\\
&=\int_{B(o,r+\varepsilon R)} (h-c)^2 |d\chi|^2\,,\end{split}\end{equation*} and 
hence
\begin{equation*}\begin{split}
\int_{B(o,r)} |dh|^2&\le \int_{B(o,r+\varepsilon R)} |d(\chi(h-c))|^2=\int_{B(o,r+\varepsilon R)} (h-c)^2 |d\chi|^2\\
&\le \frac{1}{\varepsilon^2 R^2} \int_{B(o,r+\varepsilon R)\setminus B(o,r) } (h-c)^2\\
&\le \frac{1}{\varepsilon^2 R^2} \int_{B(x,\varepsilon R+\varepsilon r) } (h-c)^2.
\end{split}\end{equation*}
The hypothesis that $\varepsilon \le 1/12$ implies that the ball $B(x,\varepsilon R+\varepsilon r)$ is remote, hence if we choose 
$$c=h_{B(x,\varepsilon (R+r)) }=\frac{1}{\volu B(x,,\varepsilon (R+r)) }\int_{B(x,,\varepsilon (R+r)) } h$$ then the Poincar\'e inequality and the fact that $r+R\le 3 R$ imply :
$$\int_{B(o,r)} |dh|^2\le 9 \mu  \int_{B(x,6\varepsilon R) } |dh|^2.$$
But we have :
$$B(x,6\varepsilon R)\subset B(o,r+6\varepsilon R)\setminus B(o,r-6\varepsilon R)\,,$$ hence
we get 
$$\int_{B(o,r-6\varepsilon R)} |dh|^2\le 9\mu  \int_{B(o,r+6\varepsilon R)\setminus B(o,r-6\varepsilon R) } |dh|^2.$$
And for all $r\in [R,R-12\varepsilon R]$ we get :
$$\int_{B(o,r)} |dh|^2\le \frac{9\mu}{1+9\mu} \int_{B(o,r+12\varepsilon R)} |dh|^2.$$
We iterate this inequality and get 
$$\int_{B(o,R)} |dh|^2\le \left(\frac{9\mu}{1+9\mu} \right)^N\int_{B(o,2 R)} |dh|^2$$
provide that 
$N 12\varepsilon R\le R$ ; hence the result with 
$C=1+\frac{1}{9\mu}$ and 
$$\kappa =\left(\frac{9\mu}{1+9\mu} \right)^{\frac{1}{12 }}\,\,.$$
\endproof
\subsection{Harmonic function with polynomial growth}
We can now prove the following extension of Kasue's results :
\begin{thm}Let $(M,g)$ be a complete Riemannian manifold whose all remote balls $B=B(x,r)$ satisfy a scale Poincar\'e inequality :
$$\forall \varphi\in \cC^1(2B)\,:\, \|\varphi-\varphi_B\|^2_{L^2(B)}\le \mu r^2(B) \|d\varphi\|^2_{L^2(2B)}$$
Assume that balls anchored at $o$ have polynomial growth :
$$\volu B(o,R)\le C R^\mu$$
 and that geodesic spheres have sublinear diameter growth :
 $$\lim_{t\to+\infty} \frac{\rho(t)}{t}=0$$
 then any harmonic function on $(M,g)$ with polynomial growth is constant.
\end{thm}
\proof Let $h\colon M\rightarrow \R$ be a harmonic function with polynomial growth :
$$h(x)\le C (1+d(o,x))^\nu .$$
We will defined
$$E_R=\int_{B(o,R)} |dh|^2\ \  \mathrm{and}
\ \ \epsilon(r)=\sup_{t\ge r} \frac{\rho(t)}{t}\,.$$
We remark first that using the cut off function $\xi$ defined by 
$$\xi(x)=\begin{cases}
1&\text{on}\,\, B(o,R)\\
\frac{2R-d(o,x)}{ R}&\text{on}\,\, B(o,2R)\setminus B(o,R)\\
0&\text{outside}\,\, B(o,2R)\\
\end{cases}$$
We obtain
\begin{equation}\label{cappo}
E_R\le \int_{B(o,2R)} |d(\xi h)|^2= \int_{B(o,2R)} |h|^2|d\xi|^2\le C R^{2\nu+\mu-2}.\end{equation}

If we iterate the inequality obtained in \lref{inefunda}, we get for all $R$ such that $\epsilon(R)\le 1/12$ :
$$E_R\le C^\ell \kappa^{\sum_{j=0}^{\ell-1} \frac{1}{\epsilon\left(2^j R\right)}}\, E_{2^\ell R}\,\,.$$

Using the estimation (\ref{cappo}), we get 
\begin{equation}\label{energy}E_R\le C(R) e^{\ell 
\left(\frac{\log\kappa}{\ell}\sum_{j=0}^{\ell-1} \frac{1}{\epsilon\left(2^j R\right)}+\log(2)(2\nu+\mu-2)+\log C\right)}\,\,.
\end{equation}
But the Cesaro theorem convergence implies that :
$$\lim_{\ell\to +\infty} \frac 1\ell\sum_{j=0}^{\ell-1} \frac{1}{\epsilon\left(2^j R\right)}=+\infty$$
hence if we let $\ell\to +\infty$ in the inequality (\ref{energy}) we get $E_R=0$  and this for all sufficiently large $R$, hence $h$ is constant.
\endproof
\subsection{Extension}
A slight variation of the arguments yields the following extension:

\begin{thm}Let $(M,g)$ be a complete Riemannian manifold whose all remote balls $B=B(x,r)$ satisfies a scale Poincar\'e inequality :
$$\forall \varphi\in \cC^1(2B)\,:\, \|\varphi-\varphi_B\|^2_{L^2(B)}\le \mu r^2(B) \|d\varphi\|^2_{L^2(2B)}$$
Assume that the  geodesic spheres have sublinear diameter growth :
 $$\lim_{t\to+\infty} \frac{\rho(t)}{t}=0\ \ \mathrm{ and\ \  let }\ \ 
\epsilon(r)=\sup_{t\ge r} \frac{\rho(t)}{t}\,\,.$$
 Let $h\colon M\rightarrow \R$ be a harmonic function and assume that 
 $\displaystyle I_R=\int_{B(o,R)} h^2$ satisfy 
 $$\log I(R)=o\left(\int_1^{R/4} \frac{dt}{t\epsilon(t)} \right)$$
 then $h$  is constant.
\end{thm}
\proof Indeed, the above argumentation shows that 
if $R$ is large enough then
$$E_R\le  M(\ell,R) I(2^{\ell+1}R)\, 4^{-\ell} R^{-2}$$
where
\begin{equation*}\begin{split}
\log(M(\ell,R))&=\log\left(C^\ell \kappa^{\sum_{j=0}^{\ell-1} \frac{1}{\epsilon\left(2^j R\right)}}\right)\\
&=\ell \log C+\log\kappa\left( \sum_{j=0}^{\ell-1} \frac{1}{\epsilon\left(2^j R\right)}\right).
\end{split}\end{equation*}
But 
$$ \sum_{j=0}^{\ell-1} \frac{1}{\epsilon\left(2^j R\right)}\ge
 \frac{1}{\log 2} \sum_{j=0}^{\ell-1} \int_{2^{j-1}R}^{2^j R}\frac{dt}{t\epsilon\left(t\right)}\ge \frac{1}{\log 2}\int_{R/2}^{2^{\ell-1} R}\frac{dt}{t\epsilon\left(t\right)}  .$$
 Hence we get the inequality :
 $$\log E_R\le \log I\!\left(2^{\ell+1}R\right)-\ell\log(4)+\ell\log C+\frac{\log \kappa}{\log 2}
   \,\int_{R/2}^{2^{\ell-1} R}\frac{dt}{t\epsilon\left(t\right)}-2\log R.$$
   It is then easy to conclude.
\endproof

\section{Lipschitz regularity of harmonic functions}
We are going to prove that a Lipschitz regularity for harmonic function analogous to the the Cheng-Yau gradient inequality :
\begin{thm}
 Let $(M^n,g)$ be a complete Riemannian manifold that satisfy the doubling condition  :
 there is a constant $\upvartheta$ such that for any $x\in M$ and radius $R>0$ :
$$\hspace{1cm} \volu B(x,2R)\le  \upvartheta\volu B(x,R)$$
and assume moreover that the  Ricci curvature  satisfies a quadratic decay lower bound
$$\hspace{1cm} \,\,\, \ricci\ge - \frac{\kappa^2}{r^2(x)} g\,\,,$$
where for a fixed point $o\in M$ : $r(x):=d(o,x)$.

Assume that the diameter of geodesic sphere growth slowly
$$\diam \partial B(o,R)= \sup_{x,y\in \partial B(o,R)} d(x,y)=o(R)$$
then there is a constant $C$ such that  for any geodesic ball $B\subset M$ and any harmonic function $h\colon 3B\rightarrow \R$
$$\sup_{x\in B} |dh|^2(x)\le \frac{C}{\volu 2B}\int_{2B}|dh|^2.$$
\end{thm}

\proof According to \cite[Proposition 5.3]{CRiesz}, we need only to show that there is a constant $C$ such that  if $R>0$ and if 
$h\colon B(o,2R)\rightarrow \R$ is a harmonic function then for any 
$s\le \sigma\le R :$
\begin{equation}\label{deG}
\frac{1}{\volu B(o,s)}\int_{B(o,s)}|dh|^2\le \frac{C}{\volu B(o,\sigma)}\int_{B(o,\sigma)}|dh|^2.\end{equation}
According to the \lref{inefunda}, for all $\eta>0$, there is a $R_0>0$ such that for all $R\ge R_0$, then
$$\int_{B(o,R)}|dh|^2\le \eta \int_{B(o,2R)} |dh|^2.$$
Hence for all $R\ge R_0:$
$$\frac{1}{\volu B(o,R)}\int_{B(o,R)} |dh|^2\le \eta\,\upvartheta\,\, \frac{1}{\volu B(o,2R)} \int_{B(o,2R)}|dh|^2.$$
Choose $\eta=\upvartheta^{-1}$, then we get 
that for all $R_0\le s\le \sigma\le R :$
$$\frac{1}{\volu B(o,s)}\int_{B(o,s)} |dh|^2\le \frac{\upvartheta}{\volu B(o,\sigma)}\int_{B(o,\sigma)}|dh|^2.$$
The Ricci curvature being bounded on $B(o,3R_0)$, the Cheng and Yau gradient estimate yields a constant $B$ such that for all
$x\in B(o,R_0)$ :
$$|dh|^2(x)\le \frac{B}{\volu B(o,2R_0)}\int_{B(o,2R_0)} |dh|^2$$
Hence the estimate (\ref{deG}) holds with $C=\max \{B, \upvartheta\}$.
\endproof

\end{document}